\documentclass{article}
\usepackage{amsmath}
\usepackage{amssymb}

\begin{document}

\begin{center}
{\Large \textbf{Integral theorems for the quaternionic\\[3mm]
 $G$-monogenic mappings}}%
\begin{equation*}
\end{equation*}

\textbf{V.~S.~Shpakivskyi and T.~S.~Kuzmenko}%
\begin{equation*}
\end{equation*}
\end{center}

\textbf{Abstract.} {\small In the paper \cite{Shpakivskiy-Kuzmenko} considered a new class of quaternionic
 mappings, so-called $G$-monogenic mappings. In this paper we prove analogues
 of classical integral theorems of the holomorphic function theory: the Cauchy
 integral theorems for surface and curvilinear integrals, and
  the Cauchy integral formula for $G$-monogenic mappings.\medskip }

 \vskip2mm

\begin{equation*}
\end{equation*}

\vspace{3mm} \textbf{0. Introduction}\vspace{3mm}%
\begin{equation*}
\end{equation*}

The Cauchy integral theorem and Cauchy integral
  formula for holomorphic functions of the complex variable are
   fundamental results of the classical complex analysis.
Analogues of these results are also important tools in the quaternionic analysis.

Maybe the first quaternionic analogues  of the mentioned results for a surface in
$\mathbb{R}^3$ are obtained by G.~Moisil and N.~Theodoresco  \cite{Moisil}.
Namely, they proved some analogues of the Cauchy's theorems for a smooth surface and
for continuously differentiable functions $f$ satisfying the equality
$$\frac{\partial f}{\partial x}i+\frac{\partial f}
{\partial y}j+
\frac{\partial f}{\partial z}k=0,
$$
where $i,j,k$ are the basis quaternionic units.

R.~Fueter  \cite{Fueter} transferred the results of
paper \cite{Moisil} for a smooth surface in $\mathbb{R}^4$ and for so-called the
\textit{regular} functions $f$ which by definition satisfies the equality
$$\frac{\partial f}{\partial t}+\frac{\partial f}{\partial x}i+\frac{\partial f}
{\partial y}j+
\frac{\partial f}{\partial z}k=0.
$$
The proofs of \cite{Moisil} and \cite{Fueter} are based on the
Stokes formula.

B.~Schuler \cite{Schuler} generalized the results of Fueter by adapting  of
 Goursat's proof of the Cauchy's theorem. Due to this, he replaced the
 condition of continuity of partial
derivatives to the
 differentiability of the real-valued
components in the Stolz sense.

 A.~Sudbery \cite{Sudbery} proved the Cauchy
 theorem and Cauchy formula  under more general conditions on a function
 and a surface. He assumed as in the complex
analysis only differentiability of real-valued components of the integrand.
 He also considered the so-called differentiable surface and the rectifiable 3-chain.

O.~Herus \cite{Gerus-2011}
proved the quaternionic Cauchy theorem for a surface in
$\mathbb{R}^3$ under the same conditions on the function as in \cite{Sudbery},
but for another class of surfaces.

Some analogues of the quaternionic Cauchy theorem and Cauchy integral formula
in the theory of  $s$-\textit{regular} functions are established in \cite{Gentili}.

In this paper for quaternionic $G$-\textit{monogenic} mappings
 we prove  analogues of the Cauchy theorem for surface and curvilinear integrals
  and the Cauchy integral formula curvilinear integral.

\vskip2mm
\textbf{1. The quaternionic $G$-monogenic mappings.}

Let $\mathbb{H(C)}$ be the quaternion algebra over the field of complex numbers
$\mathbb{C}$, whose basis consists of the unit $1$ of the algebra and of the elements
$I,J,K$
 satisfying the multiplication rules:
$$I^2=J^2=K^2=-1,\,IJ=-JI=K,\,JK=-KJ=I,\,KI=-IK=J.$$

In the algebra $\mathbb{H(C)}$ there exists another basis $\{e_1,e_2,e_3,e_4\}$:
$$e_1=\frac{1}{2}(1+iI), \quad e_2=\frac{1}{2}(1-iI), \quad
e_3=\frac{1}{2}(iJ-K), \quad e_4=\frac{1}{2}(iJ+K),$$
where $i$ is the complex imaginary unit. Multiplication table in a new
basis can be represented as

\begin{equation}\label{tabl}
\begin{tabular}{c||c|c|c|c|}
$\cdot$ & $e_1$ & $e_2$ & $e_3$ & $e_4$\\ 
\hline
\hline
$e_1$ & $e_1$ & $0$ & $e_3$ & $0$\\ 
\hline
$e_2$ & $0$ & $e_2$ & $0$ & $e_4$\\ 
\hline
$e_3$ & $0$ & $e_3$ & $0$ & $e_1$\\ 
\hline
$e_4$ & $e_4$ & $0$ & $e_2$ & $0$\\ 
\hline
\end{tabular}\,\,.
\end{equation}
The unit $1$ can be decomposed as $1=e_1+e_2$.

Consider linear functionals $f_1:\mathbb{H(C)\rightarrow\mathbb{C}}$ and $f_2:\mathbb{H(C)\rightarrow\mathbb{C}}$ satisfying the equalities
$$f_1(e_1)=f_1(e_3)=1, \qquad f_1(e_2)=f_1(e_4)=0,$$
$$f_2(e_2)=f_2(e_4)=1, \qquad f_2(e_1)=f_2(e_3)=0.$$

Let us consider the vectors
$$i_1=1=e_1+e_2, \quad i_2=a_1e_1+a_2e_2, \quad i_3=b_1e_1+b_2e_2,
\quad a_k,b_k\in\mathbb{C},\,k=1,2,$$
which are a linearly independent over the field of real numbers $\mathbb{R}$.
It means that the equality
$$\alpha_1i_1+\alpha_2i_2+\alpha_3i_3=0, \qquad
\alpha_1,\alpha_2,\alpha_3\in\mathbb{R}$$ holds if and only if
 $\alpha_1=\alpha_2=\alpha_3=0$.

In the algebra $\mathbb{H(C)}$ consider the linear span $E_3:=\{\zeta=xi_1+yi_2+zi_3:x,y,z\in\mathbb{R}\}$ generated by the
vectors $i_1,i_2,i_3$ over the field $\mathbb{R}$.
Denote $f_k(E_3):=\{f_k(\zeta) : \zeta\in E_3\}$, $k=1,2$. In what follows,
we make the following essential assumption:
$f_1(E_3)=f_2(E_3)=\mathbb{C}$.
  Obviously, it holds if and only if
at least one of the numbers in each of pairs $(a_1,b_1)$ or $(a_2,b_2)$
  belongs to
 $\mathbb{C}\setminus\mathbb{R}$.

Let us introduce the notations
$$\xi_1:=f_1(\zeta)=x+ya_1+zb_1,\qquad \xi_2:=
f_2(\zeta)=x+ya_2+zb_2.$$
Now, the element $\zeta\in E_3$ can be represented
 in the form $\zeta=\xi_1e_1+\xi_2e_2$.

 A set $S\subset\mathbb{R}^3$ is associated with the set $S_\zeta:=
\{\zeta=xi_1+yi_2+zi_3:(x,y,z)\in S\}$ in $E_3$. We also note that a
 topological property of a set $S_\zeta$ in $E_3$ understand as the same
  topological property of the set $S$ in $\mathbb{R}^3$.
For example, we will say that a curve $\gamma_\zeta\subset E_3$ is homotopic to
 the zero if $\gamma\subset\mathbb{R}^3$ is homotopic to the zero, etc.

 Let $\Omega_\zeta$ be a domain in $E_3$.

A continuous mapping $\Phi:\Omega_\zeta\rightarrow\mathbb{H(C)}$ (or $\widehat{\Phi}:\Omega_\zeta\rightarrow\mathbb{H(C)}$) is \emph{right-$G$-monogenic}
\big(or resp. \emph{left-$G$-monogenic}\big) in a domain
$\Omega_\zeta\subset E_3$,
if $\Phi$ \big(or resp. $\widehat{\Phi}$\big) is differentiable in the sense
 of the Gateaux in every point of $\Omega_\zeta$, i.~e. if for every $\zeta\in
 \Omega_\zeta$ there exists an element $\Phi'(\zeta)\in\mathbb{H(C)}$ \big(or resp. $\widehat{\Phi}'(\zeta)\in\mathbb{H(C)}$\big) such that
$$\lim\limits_{\varepsilon\rightarrow 0+0}\Big(\Phi(\zeta+\varepsilon h)-\Phi(\zeta)\Big)\varepsilon^{-1}= h\Phi'(\zeta)\quad\forall\,h\in E_3$$
$$\Biggr(\text{or resp.}\,\, \lim\limits_{\varepsilon\rightarrow 0+0}
\left(\widehat{\Phi}(\zeta+\varepsilon h)-\widehat{\Phi}(\zeta)\right)
\varepsilon^{-1}= \widehat{\Phi}'(\zeta)h\quad\forall\,h\in E_3\Biggr).$$
$\Phi'(\zeta)$ is \emph{the right Gateaux derivative} in the point $\zeta$\, and $\widehat{\Phi}'(\zeta)$ is \emph{the left Gateaux derivative} in the point $\zeta$\,.

A mapping $\Phi(\zeta)$\, (or $\widehat{\Phi}(\zeta)$) of the variable $\zeta=x+yi_2+zi_3\in\Omega_\zeta$ with differentiable real--valued components
 is right-$G$-monogenic (or resp.
left-$G$-monogenic)
if and only if the following Cauchy -- Riemann conditions are satisfied
\cite{Shpakivskiy-Kuzmenko}:

\begin{equation}\label{umova-r-K-R}
\frac{\partial \Phi}{\partial y}=i_2\frac{\partial \Phi}{\partial x}\,, \qquad
 \frac{\partial \Phi}{\partial z}=i_3\frac{\partial \Phi}{\partial x}
\end{equation}
or resp.
\begin{equation}\label{umova-l-K-R}
\frac{\partial \widehat{\Phi}}{\partial y}=\frac{\partial \widehat{\Phi}}{\partial x}
i_2\,, \qquad \frac{\partial \widehat{\Phi}}{\partial z}=\frac{\partial
\widehat{\Phi}}{\partial x}i_3\,.
\end{equation}

It follows from the decomposition of the resolvent
$$(t-\zeta)^{-1}=\frac{1}{t-\xi_1}\,e_1+\frac{1}{t-\xi_2}\,e_2,\qquad
\forall\,\,t\in \mathbb{C}:\,\,t\neq \xi_1, \,\, t\neq \xi_2$$
that the points $(x,y,z)\in\mathbb{R}^3$ corresponding to the noninvertible elements $\zeta=xi_1+yi_2+zi_3$ of the algebra $\mathbb{H(C)}$ form the
straight lines in $\mathbb{R}^3$:
$$L^1: \,x+y\text{Re}\,a_1+z\text{Re}\,b_1=0, \qquad y\text{Im}\,
a_1+z\text{Im}\,b_1=0,$$
$$L^2: \,x+y\text{Re}\,a_2+z\text{Re}\,b_2=0, \qquad y\text{Im}\,
a_2+z\text{Im}\,b_2=0$$
in the three-dimensional space $\mathbb{R}^3$.

Denote by $D_k\subset\mathbb{C}$ the image of $\Omega_\zeta$ under the mapping
$f_k$,\, $k=1,2$.
A constructive description of all right- and left-$G$-monogenic mappings
 by means of holomorphic functions of the complex variable are
  obtained in the paper \cite{Shpakivskiy-Kuzmenko}.  Namely,
 proved the theorem:

Let a domain $\Omega\subset
\mathbb{R}^{3}$ is convex in the
direction of the straight lines $L^1$, $L^2$ and $f_1(E_3)=f_2(E_3)=\mathbb{C}$.
Then any
right-$G$-monogenic mapping $\Phi:\Omega_\zeta\rightarrow\mathbb{H(C)}$ can
be expressed in the form
\begin{equation}\label{Phi-r-rozklad}
\Phi(\zeta)=F_1(\xi_1)e_1+F_2(\xi_2)e_2+F_3(\xi_1)e_3+F_4(\xi_2)e_4,
\end{equation}
$$\forall\,\zeta=xi_1+yi_2+zi_3\in\Omega_\zeta,$$ where
 $F_1, F_4$ are the certain holomorphic in
a domain $D_1$ functions of the variable $\xi_1:=x+ya_1+zb_1$, and
  $F_2, F_3$ are the certain holomorphic in
a domain $D_2$ functions of the variable $\xi_2:=x+ya_2+zb_2$.

Under the same assumptions, any
left-$G$-monogenic mapping
$\widehat{\Phi}:\Omega_\zeta\rightarrow\mathbb{H(C)}$ can
be expressed in the form
\begin{equation}\label{Phi-l-rozklad}
\widehat{\Phi}(\zeta)=F_1(\xi_1)e_1+F_2(\xi_2)e_2+F_3(\xi_2)e_3+F_4(\xi_1)e_4,
\end{equation}
where $F_n,\,n=1,2,3,4$ are defined similarly to above.

\vskip2mm
\textbf{2. Cauchy integral theorem for a curvilinear integral.}
Let $\gamma$ be a Jordan rectifiable curve in $\mathbb{R}^3$. For a continuous
mapping $\Psi:\gamma_\zeta\rightarrow\mathbb{H(C)}$
 of the form
\begin{equation}\label{Phi-form}
\Psi(\zeta)=\sum\limits_{k=1}^{4}{U_k(x,y,z)e_k}+i\sum
\limits_{k=1}^{4}{V_k(x,y,z)e_k},
\end{equation}
where $(x,y,z)\in\gamma$ and $U_k:\gamma\rightarrow\mathbb{R}$,
$V_k:\gamma\rightarrow\mathbb{R}$,
we define integrals along a Jordan rectifiable curve $\gamma_\zeta$ by
 the equalities:
$$\int\limits_{\gamma_\zeta}{d\zeta\Psi(\zeta)}:=\sum\limits_{k=1}^{4}{e_k\int
\limits_{\gamma}U_k(x,y,z)dx}+\sum\limits_{k=1}^{4}{i_2e_k\int\limits_{\gamma}
U_k(x,y,z)dy}+$$
$$+\sum\limits_{k=1}^{4}{i_3e_k\int\limits_{\gamma}U_k(x,y,z)dz}+i\sum
\limits_{k=1}^{4}{e_k\int\limits_{\gamma}V_k(x,y,z)dx}+$$
$$+i\sum\limits_{k=1}^{4}{i_2e_k\int\limits_{\gamma}V_k(x,y,z)dy}+i\sum
\limits_{k=1}^{4}{i_3e_k\int\limits_{\gamma}V_k(x,y,z)dz}$$

\noindent and

$$\int\limits_{\gamma_\zeta}{\Psi(\zeta)d\zeta}:=\sum\limits_{k=1}^{4}
{e_k\int\limits_{\gamma}U_k(x,y,z)dx}+\sum\limits_{k=1}^{4}{e_ki_2\int
\limits_{\gamma}U_k(x,y,z)dy}+$$
$$+\sum\limits_{k=1}^{4}{e_ki_3\int\limits_{\gamma}U_k(x,y,z)dz}+i\sum
\limits_{k=1}^{4}{e_k\int\limits_{\gamma}V_k(x,y,z)dx}+$$
$$+i\sum\limits_{k=1}^{4}{e_ki_2\int\limits_{\gamma}V_k(x,y,z)dy}+i\sum
\limits_{k=1}^{4}{e_ki_3\int\limits_{\gamma}V_k(x,y,z)dz},$$
where $d\zeta:=dx+i_2dy+i_3dz$.

Let $\Sigma$ be a piece-smooth surface in $\mathbb{R}^{3}$.
For a continuous
function $\Psi:\Sigma_{\zeta}\rightarrow\mathbb{H(C)}$ of the
form (\ref{Phi-form}), where $(x,y,z)\in\Sigma$, we define surface
integrals on $\Sigma_{\zeta}$ with the differential form
$\sigma:=dydz+dzdxi_{2}+dxdyi_{3}$ by the
equalities
 $$\int\limits_{\Sigma_\zeta}\hspace{-1mm}\sigma\Psi(\zeta):=
 \sum\limits_{k=1}^{4}e_{k}\hspace{-1mm}\int\limits_{\Sigma}\hspace{-1mm}
 U_{k}(x,y,z)dydz+
\sum\limits_{k=1}^{4}i_{2}e_{k}\hspace{-1mm}\int\limits_{\Sigma}\hspace{-1mm}
U_{k}(x,y,z)dzdx+$$
\vspace{-0.5mm}
$$+\sum\limits_{k=1}^{4}i_{3}e_{k}\int\limits_{\Sigma}U_{k}(x,y,z)dxdy
+i\sum\limits_{k=1}^{4}e_{k}\int\limits_{\Sigma}V_{k}(x,y,z)dydz+$$
\vspace{-0.5mm}
$$+i\sum\limits_{k=1}^{4}i_{2}e_{k}\int\limits_{\Sigma}V_{k}(x,y,z)dzdx+
i\sum\limits_{k=1}^{4}i_{3}e_{k}\int\limits_{\Sigma}V_{k}(x,y,z)dxdy;$$

$$\int\limits_{\Sigma_\zeta}\hspace{-1mm}\Psi(\zeta)\sigma:=
 \sum\limits_{k=1}^{4}e_{k}\hspace{-1mm}\int\limits_{\Sigma}\hspace{-1mm}
 U_{k}(x,y,z)dydz+
\sum\limits_{k=1}^{4}e_{k}i_{2}\hspace{-1mm}\int\limits_{\Sigma}\hspace{-1mm}
U_{k}(x,y,z)dzdx+$$
\vspace{-0.5mm}
$$+\sum\limits_{k=1}^{4}e_{k}i_{3}\int\limits_{\Sigma}U_{k}(x,y,z)dxdy
+i\sum\limits_{k=1}^{4}e_{k}\int\limits_{\Sigma}V_{k}(x,y,z)dydz+$$
\vspace{-0.5mm}
$$+i\sum\limits_{k=1}^{4}e_{k}i_{2}\int\limits_{\Sigma}V_{k}(x,y,z)dzdx+
i\sum\limits_{k=1}^{4}e_{k}i_{3}\int\limits_{\Sigma}V_{k}(x,y,z)dxdy.$$

If a function $\Psi:\Omega_{\zeta}\rightarrow\mathbb{H(C)}$ is
continuous together with partial derivatives of the first order in
a domain $\Omega_{\zeta}$, and $\Sigma$ is a piece-smooth surface
in $\Omega$, and the edge $\gamma$ of surface $\Sigma$ is a
rectifiable Jordan curve, then the following analogues of the
Stokes formula are true:
$$\int\limits_{\gamma_{\zeta}}d\zeta\Psi(\zeta)=\int\limits_{\Sigma_{\zeta}}
\left(i_2\frac{\partial\Psi}{\partial
x}-\frac{\partial\Psi}{\partial y}\right)dxdy
+\left(i_3\frac{\partial\Psi}{\partial
y}-i_2\frac{\partial\Psi}{\partial z}\right)dydz+$$
\begin{equation}\label{form-Stoksa}+\left(\frac{\partial\Psi}{\partial
z}-i_3\frac{\partial\Psi}{\partial x}\right)dzdx,
\end{equation}

$$\int\limits_{\gamma_{\zeta}}\Psi(\zeta)d\zeta=\int\limits_{\Sigma_{\zeta}}
\left(\frac{\partial\Psi}{\partial
x}i_{2}-\frac{\partial\Psi}{\partial y}\right)dxdy
+\left(\frac{\partial\Psi}{\partial
y}i_{3}-\frac{\partial\Psi}{\partial z}i_{2}\right)dydz+$$
\begin{equation}\label{form-Stoksa-1}+\left(\frac{\partial\Psi}{\partial
z}-\frac{\partial\Psi}{\partial x}i_{3}\right)dzdx.
\end{equation}

Now, the next theorem is a result of the formulae
(\ref{form-Stoksa}), (\ref{form-Stoksa-1}) and the equalities (\ref{umova-r-K-R}), (\ref{umova-l-K-R}), respectively.

\textbf{Theorem 1.} \emph{Suppose that $\Phi:\Omega_\zeta\rightarrow
\mathbb{H(C)}$ is a right-$G$-monogenic mapping in a domain $\Omega_\zeta$ and $\widehat{\Phi}:\Omega_\zeta\rightarrow\mathbb{H(C)}$
is a left-$G$-monogenic mapping in $\Omega_\zeta$\,.
Suppose also that $\Sigma$ is a
piece-smooth surface in $\Omega$, and the edge $\gamma$ of surface
$\Sigma$ is a rectifiable Jordan curve. Then}
\begin{equation}\label{form-Koshi-po-kryv}
\int\limits_{\gamma_\zeta}d\zeta\,\Phi(\zeta)=\int\limits_{\gamma_\zeta}
\widehat{\Phi}(\zeta)d\zeta=0.
\end{equation}

In the case where a domain $\Omega$ is convex, then by the usual way
 (see, e.~g., \cite{Privalov}) the equality (\ref{form-Koshi-po-kryv})
  can be prove for an arbitrary closed Jordan rectifiable curve $\gamma_\zeta$.

In the case where a domain $\Omega$ is an arbitrary, then similarly
 to the proof of Theorem 3.2 \cite{Blum} we can
prove the following

\vskip 1mm
\textbf{Theorem 2.} \emph{Let $\Phi:\Omega_\zeta\rightarrow\mathbb{H(C)}$ be a
right-$G$-monogenic mapping in a domain $\Omega_\zeta$ and $\widehat{\Phi}:\Omega_\zeta\rightarrow\mathbb{H(C)}$
be a left-$G$-monogenic mapping in $\Omega_\zeta$\,. Then for every
closed Jordan rectifiable curve $\gamma_\zeta$ homotopic to a point in
$\Omega_\zeta$\,, the the equalities \em (\ref{form-Koshi-po-kryv}) \em holds.}

\vskip2mm
\textbf{3. Cauchy integral formula.} To establish the Cauchy integral formula
for a curvilinear integral, consider
the following auxiliary statement:

\vskip 1mm
\textbf{Lemma.}  \emph{Suppose that a domain $\Omega\subset
\mathbb{R}^3$ is convex in the direction of the straight lines $L^1,\,L^2$ and $f_1(E_3)=f_2(E_3)=\mathbb{C}$. Suppose also that $\Phi:\Omega_\zeta\rightarrow
\mathbb{H(C)}$ is a right-$G$-monogenic mapping in $\Omega_\zeta$, and $\widehat{\Phi}:\Omega_\zeta\rightarrow\mathbb{H(C)}$ is a left-$G$-monogenic
 mapping in $\Omega_\zeta$, and $\gamma_\zeta$ is an arbitrary rectifiable curve
 in $\Omega_\zeta$.
 Then}
\begin{equation}\label{lema-1-r}
\int\limits_{\gamma_\zeta}d\zeta\,\Phi(\zeta)=e_1\int\limits_{\gamma_1}F_1(\xi_1)
d\xi_1+e_2\int\limits_{\gamma_2}F_2(\xi_2)d\xi_2+e_3\int\limits_{\gamma_1}
F_3(\xi_1)d\xi_1+e_4\int\limits_{\gamma_2}F_4(\xi_2)d\xi_2\,,
\end{equation}

\noindent\emph{and respectively}
\begin{equation}\label{lema-1-l}
\int\limits_{\gamma_\zeta}\widehat{\Phi}(\zeta)d\zeta=e_1\int\limits_{\gamma_1}
F_1(\xi_1)d\xi_1+e_2\int\limits_{\gamma_2}F_2(\xi_2)d\xi_2+e_3\int\limits_{\gamma_2}
F_3(\xi_2)d\xi_2+e_4\int\limits_{\gamma_1}F_4(\xi_1)d\xi_1\,,
\end{equation}
\emph{where $\gamma_k$ is the image of $\gamma_\zeta$ under the mapping $f_k$
and  $F_n$ is the same function as in}
(\ref{Phi-r-rozklad}) \emph{and} (\ref{Phi-l-rozklad}) \emph{respectively}.
\vskip 1mm

\textbf{Proof.} The equality (\ref{lema-1-r}) follows immediately from the
representation (\ref{Phi-r-rozklad}), the equality $d\zeta=d\xi_1e_1+d\xi_2e_2$
 and the multiplication rules (\ref{tabl}). Similarly we can prove the equality
  (\ref{lema-1-l}). Lemma is proved.

Let $\zeta\in E_3$. An inverse
  element $\zeta^{-1}$ is of the following form:
\begin{equation}\label{zeta-1}
\zeta^{-1}=\frac{1}{\xi_1}e_1+\frac{1}{\xi_2}e_2
\end{equation}
and it exists if and only if $\xi_1\neq 0$ and  $\xi_2\neq 0$.

Let $\zeta_0=\xi^{(0)}_{1}e_1+\xi^{(0)}_{2}e_2$ be a point in a domain
$\Omega_\zeta\subset E_3$. In a neighborhood of $\zeta_0$ contained in
 $\Omega_\zeta$ let us take a circle $C(\zeta_0)$ with the center at the point
  $\zeta_0$. By $C_k\subset\mathbb{C}$ we denote the image of $C(\zeta_0)$
  under the mapping $f_k$,\, $k=1,2$.  We assume that
   the circle $C(\zeta_0)$ \emph{embraces the set} $\{\zeta-\zeta_0:
   \zeta\in L^1_\zeta\cup L^2_\zeta\}$. It means that $C_k$ bounds some
    domain $D_k'$ and
    $\xi^{(0)}_{k}\in D_k'$, \, $k=1,2$.

We say that the curve $\gamma_\zeta\subset\Omega_\zeta$ \emph{embraces once the set} $\{\zeta-\zeta_0:\zeta\in L^1_\zeta\cup L^2_\zeta\}$, if there exists a circle
$C(\zeta_0)$
which embraces the mentioned set and is homotopic to $\gamma_\zeta$ in the domain $\Omega_\zeta\setminus\{\zeta-\zeta_0:\zeta\in L^1_\zeta\cup L^2_\zeta\}$.

\vskip 1mm
\textbf{Theorem 3.} \emph{Suppose that a domain $\Omega \subset \mathbb{R}^3$
 is
convex in the direction of the straight lines $L^1,\,L^2$ and $f_1(E_3)=f_2(E_3)=
\mathbb{C}$. Suppose also that $\Phi:\Omega_\zeta\rightarrow\mathbb{H(C)}$ is a
 right-$G$-monogenic mapping in $\Omega_\zeta$ and $\widehat{\Phi}:
 \Omega_\zeta\rightarrow
 \mathbb{H(C)}$ is a left-$G$-monogenic mapping in $\Omega_\zeta$. Then for
 every point $\zeta_0\in\Omega_\zeta$ the following equalities are true:}
\begin{equation}\label{r-int-formula}
\Phi(\zeta_0)=\frac{1}{2\pi i}\int\limits_{\gamma_\zeta}(\zeta-\zeta_0)^{-1}
d\zeta\,\Phi(\zeta)
\end{equation}
\emph{and}
\begin{equation}\label{l-int-formula}
\widehat{\Phi}(\zeta_0)=\frac{1}{2\pi i}\int\limits_{\gamma_\zeta}\widehat{\Phi}(\zeta)(\zeta-\zeta_0)^{-1}d\zeta,
\end{equation}
\emph{where $\gamma_\zeta$ is an arbitrary closed Jordan rectifiable curve in
 $\Omega_\zeta$, that embraces once the set $\{\zeta-\zeta_0:\zeta\in L_1\cup L_2\}$.}
\vskip 1mm

\textbf{Proof.} Inasmuch as $\gamma_\zeta$ is homotopic to $C(\zeta_0)$ in the domain $\Omega_\zeta\setminus\{\zeta-\zeta_0:\zeta\in L_1\cup L_2\}$, it follows
from Theorem 2 that
$$\frac{1}{2\pi i}\int\limits_{\gamma_\zeta}(\zeta-\zeta_0)^{-1}d\zeta\,\Phi(\zeta)
=\frac{1}{2\pi i}\int\limits_{C(\zeta_0)}(\zeta-\zeta_0)^{-1}d\zeta\,\Phi(\zeta).$$

Further, using the equality (\ref{zeta-1}), Lemma
and the integral Cauchy formula for holomorphic functions $F_n$, we obtain immediately
 the following equalities:
$$\frac{1}{2\pi i}\int\limits_{C(\zeta_0)}(\zeta-\zeta_0)^{-1}d\zeta\,\Phi(\zeta)
=e_1\frac{1}{2\pi i}\int\limits_{C_1}\frac{F_1(\xi_1)}{\xi_1-\xi_1^{(0)}}d\xi_1+
e_2\frac{1}{2\pi i}\int\limits_{C_2}\frac{F_2(\xi_2)}{\xi_2-\xi_{2}^{(0)}}d\xi_2+$$
$$+e_3\frac{1}{2\pi i}\int\limits_{C_1}\frac{F_3(\xi_1)}{\xi_1-\xi^{(0)}_{1}}d\xi_1+e_4\frac{1}{2\pi i}\int\limits_{C_2}\frac{F_4(\xi_2)}{\xi_2-\xi^{(0)}_{2}}d\xi_2=$$
$$=F_1\big(\xi^{(0)}_{1}\big)e_1+F_2\big(\xi^{(0)}_{2}\big)e_2+F_3\big(\xi^{(0)}_{1}
\big)e_3+F_4\big(\xi^{(0)}_{2}\big)e_4=\Phi(\zeta_0),$$
where $\zeta_0=\xi^{(0)}_{1}e_1+\xi^{(0)}_{2}e_2$.
Similarly can be proved the equality (\ref{l-int-formula}). The theorem is proved.

We note that the method of this proof is similarly to the proof
of Theorem 6 of the paper \cite{Plaksa_Puht}, where Cauchy integral formula
is obtained in a finite-dimensional semi-simple commutative algebra.

\vskip2mm
\textbf{4. Cauchy integral theorem for a surface integral.}
Let $\Omega$ be a bounded domain in $\mathbb{R}^3$. For a continuous mapping $\Psi:\Omega_\zeta\rightarrow\mathbb{H(C)}$  of the form (\ref{Phi-form}),
where $(x,y,z)\in\Omega$ and $U_k:\Omega\rightarrow\mathbb{R},\,V_k:
\Omega\rightarrow\mathbb{R}$, we define a volume integral by the equality
$$\int\limits_{\Omega_\zeta}{\Phi(\zeta)dxdydz}:=\sum\limits_{k=1}^{4}
{e_k\int\limits_{\Omega}U_k(x,y,z)dxdydz}+i\sum\limits_{k=1}^{4}
{e_k\int\limits_{\Omega}V_k(x,y,z)dxdydz}.$$

Let $\Sigma$ be a piece-smooth surface in $\mathbb{R}^3$. For a
 continuous mapping $\Psi:\Sigma_\zeta\rightarrow\mathbb{H(C)}$  of the form
  (\ref{Phi-form}),
where $(x,y,z)\in\Sigma$ and $U_k:\Sigma\rightarrow\mathbb{R},\,V_k:
\Sigma\rightarrow\mathbb{R}$, we
 define the surface integrals on a piece-smooth surface $\Sigma_\zeta$ with
 the differential form $\sigma:=dydz+dzdxi_2+dxdyi_3$ by the equalities
$$\int\limits_{\Sigma_\zeta}{\sigma\Psi(\zeta)}:=\sum\limits_{k=1}^{4}
{e_k\int\limits_{\Sigma}U_k(x,y,z)dydz}+\sum\limits_{k=1}^{4}{i_2e_k\int
\limits_{\Sigma}U_k(x,y,z)dzdx}+$$
$$+\sum\limits_{k=1}^{4}{i_3e_k\int\limits_{\Sigma}U_k(x,y,z)dxdy}+i\sum
\limits_{k=1}^{4}{e_k\int\limits_{\Sigma}V_k(x,y,z)dydz}+$$
$$+i\sum\limits_{k=1}^{4}{i_2e_k\int\limits_{\Sigma}V_k(x,y,z)dzdx}+i\sum
\limits_{k=1}^{4}{i_3e_k\int\limits_{\Sigma}V_k(x,y,z)dxdy}$$
and
$$\int\limits_{\Sigma_\zeta}{\Psi(\zeta)\sigma}:=\sum\limits_{k=1}^{4}
{e_k\int\limits_{\Sigma}U_k(x,y,z)dydz}+\sum\limits_{k=1}^{4}{e_ki_2\int
\limits_{\Sigma}U_k(x,y,z)dzdx}+$$
$$+\sum\limits_{k=1}^{4}{e_ki_3\int\limits_{\Sigma}U_k(x,y,z)dxdy}+i\sum
\limits_{k=1}^{4}{e_k\int\limits_{\Sigma}V_k(x,y,z)dydz}+$$
$$+i\sum\limits_{k=1}^{4}{e_ki_2\int\limits_{\Sigma}V_k(x,y,z)dzdx}+i\sum
\limits_{k=1}^{4}{e_ki_3\int\limits_{\Sigma}V_k(x,y,z)dxdy}.$$

If a domain $\Omega\subset \mathbb{R}^3$ has a closed piece-smooth boundary
$\partial\Omega$ and a mapping $\Psi:\Omega_\zeta\rightarrow\mathbb{H(C)}$
 is continuous together with partial derivatives of the first order up to the boundary $\partial\Omega_\zeta$, then the following analogues of the Gauss -- Ostrogradsky
  formula are true:
\begin{equation}\label{Gauss-Ostr-r}
\int\limits_{\partial\Omega_\zeta}{\sigma\Psi(\zeta)}=\int\limits_{\Omega_\zeta}
\left(\frac{\partial\Psi}{\partial x}+i_2\frac{\partial\Psi}{\partial y}+i_3\frac{\partial\Psi}{\partial z}\right)dxdydz.
\end{equation}
and
\begin{equation}\label{Gauss-Ostr-l}
\int\limits_{\partial\Omega_\zeta}{\Psi(\zeta)\sigma}=\int
\limits_{\Omega_\zeta}\left(\frac{\partial\Psi}{\partial x}+\frac{\partial\Psi}
{\partial y}i_2+\frac{\partial\Psi}{\partial z}i_3\right)dxdydz.
\end{equation}

Now, the next theorem is a result of the formulas (\ref{Gauss-Ostr-r}),
 (\ref{Gauss-Ostr-l}) and the conditions (\ref{umova-r-K-R}), (\ref{umova-l-K-R}),
 respectively.

 \vskip2mm
\textbf{Theorem 4.} \emph{Suppose that $\Omega$ has a
closed piece-smooth boundary $\partial\Omega$. Suppose also that the mapping $\Phi:\Omega_\zeta\rightarrow\mathbb{H(C)}$ is a right-$G$-monogenic in $\Omega_\zeta$
 and $\widehat{\Phi}:\Omega_\zeta\rightarrow\mathbb{H(C)}$ is a left-$G$-monogenic
in $\Omega_\zeta$, and these mappings are continuous together with partial
derivatives of the first order up to the boundary $\partial\Omega_\zeta$. Then}
\begin{equation}\label{cauchy-surface-r}
\int\limits_{\partial\Omega_\zeta}\sigma\Phi(\zeta)=\int\limits_{\Omega_\zeta}
(1+i_2^2+i_3^2)\Phi'(\zeta)dxdydz
\end{equation}
\emph{and}
\begin{equation}\label{cauchy-surface-l}
\int\limits_{\partial\Omega_\zeta}\widehat{\Phi}(\zeta)\sigma=\int
\limits_{\Omega_\zeta}\widehat{\Phi}'(\zeta)(1+i_2^2+i_3^2)dxdydz.
\end{equation}

\vskip2mm
\textbf{Corollary.} \emph{Under the conditions of Theorem
\em 4
 \em
  with the additional assumption $1+i_2^2+i_3^2=0$, i.~e. the mappings $\Phi$ and
 $\widehat{\Phi}$  are
 solutions of the three-dimensional Laplace equation, then the equalities}
  (\ref{cauchy-surface-r})
 \emph{and} (\ref{cauchy-surface-l}) \emph{can be rewritten in the form}
$$\int\limits_{\partial\Omega_\zeta}\sigma\Phi(\zeta)=\int\limits_{\partial
\Omega_\zeta}\widehat{\Phi}(\zeta)\sigma=0.$$

\vskip7mm
Vitalii Shpakivskyi

 Department of Complex Analysis and Potential Theory

 Institute of Mathematics of the National Academy of Sciences of
Ukraine,

 3, Tereshchenkivs'ka st.

01601 Kyiv-4

 UKRAINE

 http://www.imath.kiev.ua/\~{}complex/

 \ e-mail: shpakivskyi@mail.ru,\, shpakivskyi@imath.kiev.ua

\vspace{8mm}
Tetyana Kuzmenko

Faculty of mathematics and physics

Zhytomyr State Ivan Franko University

40, Velyka Berdychivs'ka st.

10008 Zhytomyr

  UKRAINE

\ e-mail: kuzmenko.ts@mail.ru

\end{document}